\documentclass[12pt,a4paper]{amsart}

\usepackage{amsmath}
\usepackage{amssymb}
\usepackage{amsthm}
\usepackage{latexsym}
\usepackage[T1]{fontenc}
\usepackage{pgf,tikz}
\usetikzlibrary{decorations.pathreplacing}
\usetikzlibrary{arrows}
\usepackage{float}
\usepackage{algorithm}
\usepackage[noend]{algpseudocode}

\theoremstyle{plain}
\newtheorem{theorem}{Theorem}[section]
\newtheorem{corollary}[theorem]{Corollary}
\newtheorem{lemma}[theorem]{Lemma}
\newtheorem{proposition}[theorem]{Proposition}

\theoremstyle{definition}
\newtheorem{definition}[theorem]{Definition}%[section]

\theoremstyle{remark}
\newtheorem*{remark}{Remark}
\newtheorem{example}{Example}[section]

\numberwithin{equation}{section}

\numberwithin{table}{section}
\def\lbr{\underset{\textup{\textsf L}}{\sim}}
\def\sbr{\underset{\textup{\textsf S}}{\sim}}

\def\ds{\displaystyle}
\def\ss{{\mathcal S}}

\def\Red{\text{R}}
\def\sup{\text{sup}}
\def\tt{\text{T}}
\def\Des{\text{Des}}
\def\wo{w_{\text{o}}}
\def\SBT{\text{SBT}}
\def\bb{\mathfrak{b}}
\def\cc{\mathfrak{c}}

\begin{document}

\title[]{Diameter of the commutation classes graph of a permutation}

\author{Gon\c{c}alo Gutierres}
\address{University of Coimbra, CMUC, Department of Mathematics,
ORCiD 0000-0001-9480-498X} \email{ggutc@mat.uc.pt}

\author{Ricardo Mamede}
\address{University of Coimbra, CMUC, Department of Mathematics,
ORCiD 0000-0002-9264-6604} \email{mamede@mat.uc.pt}

\author{Jos\'{e} Luis Santos}
\address{University of Coimbra, CMUC, Department of Mathematics,
ORCiD 0000-0002-2727-6774} \email{zeluis@mat.uc.pt}

\thanks{This work was partially supported by the Centre for Mathematics of the University of Coimbra - UIDB/00324/2020, funded by the Portuguese Government through FCT/MCTES}


%\subjclass[2000]{05A17, 05E05, 05E10, 68Q17}

\maketitle
\begin{abstract}

We define a statistic on the graph of commutation classes of a permutation of the symmetric group which is used to show that these graphs are equipped with a ranked poset structure, with a minimum and maximum. This characterization also allows us to compute the diameter of the commutation graph for any permutation, from which the results for the longest permutation and for fully commutative permutations are recovered.
\end{abstract}

\section{Introduction}

Given an integer $n\geq 2$, we let $\ss_n$ denote the symmetric group on the alphabet $[n]=\{1,\ldots,n\}$, with composition of permutations performed from right to left. We  usually write permutations in one-line notation $w=w_1w_2\cdots w_n$, where $w_i=w(i)$.

The symmetric group $\ss_n$ is an example of the more general concept of a Coxeter group \cite{brenti,coxeter}, which are groups  $G$ that can be generated by a set $S=\{s_1,\ldots,s_m\}\subset G$ satisfying relations $(s_is_j)^{m_{ij}}=1$, where $m_{ii}=1$ and $m_{ij}\geq 2$ for $i\neq j$. Any element $w\in G$ can be written as a finite product of elements of $S$. If $w=s_{i_1}s_{i_2}\cdots s_{i_{\ell}}$ with $\ell$ minimal, the word $i_1i_2\cdots i_{\ell}$ is called a {\it reduced word} (or {\it reduced decomposition}) of $g$. In this case, we define the {\it length} of $w$ by $\ell(w)=\ell$.     The set of all reduced words of $w$ is denoted by $\Red(w)$.

The symmetric group has a Coxeter representation with generators $s_i$, the adjacent transposition interchanging the elements $i$ and $i+1$, for $1\leq i\leq n-1$, which satisfy the Coxeter relations
\begin{align}
 &s_is_j=s_js_i\text{ for }|i-j|\geq 2,\label{eq:def1}\\
 &s_is_{i+1}s_1=s_{i+1}s_is_{i+1}\text{ for }1\leq i\leq n-2,\label{eq:def2}
\end{align}
and $s_i^2=1$, the identity element. The relations \eqref{eq:def1} are known as {\it commutations} or {\it short braid relations}, and the relations \eqref{eq:def2} are called {\it long braid relations}.

The graph $G(w)$, having vertex set $\Red(w)$ and an edge connecting  two reduced words if they
differ by a single Coxeter relation, has been considered by several authors. A theorem of Tits \cite{tits}  shows that $G(w)$ is connected, and  Elnitsky \cite{elnitsky}  proved that it is bipartite (see also \cite{bergeron}).  The case of the longest permutation $\wo=n\cdots 21$ has been particularly well studied. The number of reduced words in $\Red(\wo)$ was first computed algebraically by Stanley \cite{stanley} using generating functions, and later proved bijectively by Edelman and Greene \cite{balanced1}, by establishing a bijection between reduced words for $\wo$ and balanced labelings of the Rothe diagram of $\wo$. Fomin, Greene, Reiner and Shimozono \cite{reiner} generalized this result, proving a one-to-one correspondence between reduced words for $w\in\ss_n$ and standard balanced labellings of the Rothe diagram of $w$.

The diameter of $G(w)$ was first studied asymptotically by Autord and Dehornoy \cite{Dehornoy}, who showed that the diameter grows asymptotically in  ${\mathcal O}(n^4)$,  and then exactly by  Reiner and  Roichman \cite{reiner_roichman} for the longest permutation $\wo$, using hyperplane arrangements. This result was later obtained by Assaf \cite{sami} using balanced tableaux.

Contracting the commutation edges of $G(w)$ leads to the associated graph $C(w)$, known as the commutation graph of $w$, which has also received some attention.  This graph has been studied in the context of the {\it higher Bruhat order} $B(n,2)$  \cite{felsner,manin,ziegler}, and in connection with rhombic tilings of certain polygons  \cite{elnitsky}. The diameter of $C(\wo)$ has been computed in \cite{ziegler}, and a connection to geometric representation theory was explored in \cite{kim}.
 Fishel {\it et al.} \cite{fishel} gave upper and lower bounds for the cardinality of $\Red(w)$
in terms of the number of vertices in $C(w)$ and $B(w)$, where $B(w)$ is the graph obtained from $G(w)$ by contracting the long braid edges. Quotients of $\Red(w)$ by more general Coxeter relations were considered by Bergeron {\it et al.} \cite{bergeron}.

In this paper, we establish a statistic on the classes of $C(w)$, inducing a rank poset structure on $C(w)$ with a unique minimal and a unique maximal element.   This allows us to give a precise formula for the  diameter of the graph $C(w)$. We recover, as special cases, the diameter of the commutation graph for the longest element $\wo$ and the characterization of fully commutative permutations obtained by Billey, Jockusch and Stanley \cite{321Avoiding}.

\section{Reduced Words and Balanced Tableaux}

The length of a permutation $w\in\ss_n$ can also be given by the number of {\it inversions} of $w$ \cite{brenti}, that is the number of pairs $(w_j,w_i)$ such that $i<j$ and $w_j<w_i$:
$$\ell(w)=\big|\{ (w_j,w_i): i<j\text{ and }w_j<w_i\}\big|.$$
Thus, if $w=s_{i_1}s_{i_2}\cdots s_{i_{\ell(w)}}$, the word $i_1i_2\cdots i_{\ell(w)}$ is a reduced word for $w$. We consider adjacent transpositions $s_{i_j}$ acting on positions $i_j$ and $i_j+1$, and perform their composition left to right to mirror the composition of permutations. For example,  the permutation $s_1s_2=231$ acts on $123$ as
\begin{align*}
 s_1s_2&=123\cdot s_1s_2\\
 &=213\cdot s_2\\
&= 231.
\end{align*}

Another useful notion is the {\it descent set} of a permutation, which is defined as the set
$$\Des(w)=\{i:w_i>w_{i+1}\}.$$
The elements of the descent set of $w$ are called {\it descents} and can be used to obtain a reduced word for $w$. Start with $w^0=w$ and construct a sequence of permutations
\begin{equation}\label{eq:algforreducedword}
 w^0,w^1,\ldots,w^{\ell(w)},
\end{equation}
with $w^{j+1}=s_{i_j}\cdot w^j$ where $i_j$ is an element of $\Des(w^j)$, for $j\in[n]$.  Since $i_j\in \Des(w^j)$, we have $\ell(w^{j+1})=\ell(w^j)-1$, and thus the last permutation $w^{\ell(w)}$ in  \eqref{eq:algforreducedword} is the identity. Therefore, the sequence $s_{i_{\ell(w)-1}}\cdots s_{i_1}s_{i_0}$ is a reduced word for $w$.

There are, in general, several possibilities for the index $i_j\in\Des(w^j)$ in step $j$ of the procedure above. We denote by $a_{\min}$ (resp., $a_{\max}$) the reduced word for $w$ obtained by  choosing, in each step, $j$ the smallest (resp., greatest) index in $\Des(w^j)$.
As we shall see in Lemma \ref{lemma:min_max}, every braid relation in the reduced word $a_{\min}$ (resp., $a_{\max}$)  has the form $i(i+1)i$ (resp., $(i+1)i(i+1)$).

\begin{example}
 The permutation $w=25431$ has length $\ell(w)=7$, corresponding to the inversions $(2,1),(3,1),(4,1),(5,1),(4,3),(5,3)$ and $(5,4)$. We construct the reduced word $a_{\max}=4321434$ for $w$ using the procedure above, displaying the sequence of permutations in Table \ref{ex:tableex1}, where the labels on the left column (typed in italics) are descents (in this case the largest one) of the permutation in the line above. The action of each adjacent transposition $s_i$, corresponding to descent $i$, is illustrated by the numbers in bold.
 \begin{table}[htp]
\begin{center}
\begin{tabular}{c|ccccc|c}
{\it 4}&2&5&4&{\bf 3}&{\bf 1}&$w^0$\\
\cline{2-6}
{\it 3}&2&5&{\bf 4}&{\bf 1}&3&$w^1$\\
{\it 2}&2&{\bf 5}&{\bf 1}&4&3&$w^2$\\
{\it 1}&{\bf 2}&{\bf 1}&5&4&3&$w^3$\\
{\it 4}&1&2&5&{\bf 4}&{\bf 3}&$w^4$\\
{\it 3}&1&2&{\bf 5}&{\bf 3}&4&$w^5$\\
{\it 4}&1&2&3&{\bf 5}&{\bf 4}&$w^6$\\
&1&2&3&4&5&$w^7$
\end{tabular}
\end{center}
\caption{A reduced word for $25431$.}\label{ex:tableex1}
\end{table}
\end{example}

The {\it Rothe diagram} of a permutation $w\in \ss_n$, denoted by $\mathbb{D}(w)$ is  the subset of cells in the
first quadrant of the plane defined by
$$\mathbb{D}(w)=\{(i,w_j): i<j\text{ and } w_i>w_j\}\subseteq [n]\times [n].$$
The cells of $\mathbb{D}(w)$ correspond to inversions in $w$, namely
$(p,q)\in \mathbb{D}(w)$ if and only if $(w_p,q)$ is an inversion of $w$.  Therefore,
the Rothe diagram of $w$ gives a graphical representation of the inversion pairs of $w$. In particular,  the number of cells in $\mathbb{D}(w)$ is the length $\ell(w)$ of $w$.

\begin{figure}[H]
\begin{center}
 \begin{tikzpicture}
 \draw (-4,0) -- (-4,3);
 \draw (-4,0) -- (-1,0);
 \node at (-3.7, -.3) {\tiny 1};
 \node at (-3.1, -.3) {\tiny 2};
 \node at (-2.5, -.3) {\tiny 3};
 \node at (-1.9, -.3) {\tiny 4};
 \node at (-1.3, -.3) {\tiny 5};
 \node at (-4.3, .3) {\tiny 2};
 \node at (-4.3, .9) {\tiny 5};
 \node at (-4.3, 1.5) {\tiny 4};
 \node at (-4.3, 2.1) {\tiny 3};
 \node at (-4.3, 2.7) {\tiny 1};
 \draw (-4,0) rectangle (-3.4,.6);
 \draw (-4,.6) rectangle (-3.4,1.2);
 \draw (-4,1.2) rectangle (-3.4,1.8);
 \draw (-4,1.8) rectangle (-3.4,2.4);
 \draw (-2.8,.6) rectangle (-2.2,1.2);
 \draw (-2.8,1.2) rectangle (-2.2,1.8);
 \draw (-2.2,.6) rectangle (-1.6,1.2);
 %%%%
 \draw (1,0) -- (1,3);
 \draw (1,0) -- (3.9,0);
 \node at (1.3, -.3) {\tiny 1};
 \node at (1.9, -.3) {\tiny 2};
 \node at (2.5, -.3) {\tiny 3};
 \node at (3.1, -.3) {\tiny 4};
 \node at (3.6, -.3) {\tiny 5};
 \node at (0.7, .3) {\tiny 5};
 \node at (0.7, .9) {\tiny 4};
 \node at (0.7, 1.5) {\tiny 3};
 \node at (0.7, 2.1) {\tiny 2};
 \node at (0.7, 2.7) {\tiny 1};
\draw (1,0) rectangle (1.6,.6);
\draw (1,.6) rectangle (1.6,1.2);
\draw (1,1.2) rectangle (1.6,1.8);
\draw (1,1.8) rectangle (1.6,2.4);
\draw (1.6,0) rectangle (2.2,.6);
\draw (1.6,.6) rectangle (2.2,1.2);
\draw (1.6,1.2) rectangle (2.2,1.8);
\draw (2.2,0) rectangle (2.8,.6);
\draw (2.2,.6) rectangle (2.8,1.2);
\draw (2.8,0) rectangle (3.4,.6);
 \end{tikzpicture}
\end{center}
\caption{Rothe diagrams for $25431$ and $54321$.}\label{rothe1}
\end{figure}
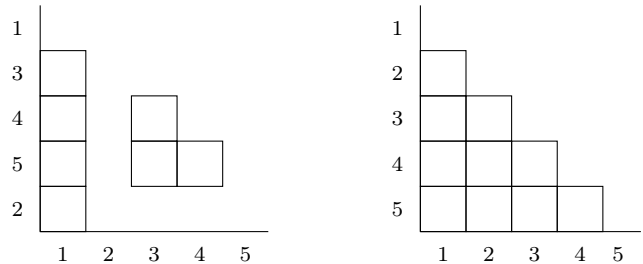

The Rothe diagram for $w=w_1w_2\cdots w_n$ can be obtained by writing $w$ vertically along the $y$-axis,
with $w_i$ at height $i$, and writing the positive numbers along the $x$-axis. Then, with this numerical arrangement of rows and columns, place a cell in
position $(w_j,i)$ whenever this is an inversion pair for $w$, for each $i,j\in{n}$. For instance, the Rothe diagrams for the permutations $25431$ and $54321$ are given in Figure \ref{rothe1}.

Note that the Rothe diagram of $w$ have the {\it transitive} property: if $(z,y)$ and $(y,x)$ are cells in $\mathbb{D}(w)$, then $(z,x)$ is also in $\mathbb{D}(w)$. These three cells form a hook with end cells $(z,y)$ and $(y,x)$.

\begin{figure}[H]
\begin{center}
 \begin{tikzpicture}
 \draw (-4,0) -- (-4,3);
 \draw (-4,0) -- (-1,0);
 \draw (-4,0) rectangle (-3.4,.6);
 \draw (-4,.6) rectangle (-3.4,1.2);
 \draw (-4,1.2) rectangle (-3.4,1.8);
 \draw (-4,1.8) rectangle (-3.4,2.4);
 \draw (-2.8,.6) rectangle (-2.2,1.2);
 \draw (-2.8,1.2) rectangle (-2.2,1.8);
 \draw (-2.2,.6) rectangle (-1.6,1.2);
 \node at (-3.7,.3) {4};
 \node at (-3.7,.9) {5};
 \node at (-3.7,1.5) {6};
 \node at (-3.7,2.1) {7};
 \node at (-2.5,.9) {2};
 \node at (-2.5,1.5) {1};
 \node at (-1.9,.9) {3};
 %%%%
  \draw (1,0) -- (1,3);
 \draw (1,0) -- (3.9,0);
 \draw (1,0) rectangle (1.6,.6);
\draw (1,.6) rectangle (1.6,1.2);
\draw (1,1.2) rectangle (1.6,1.8);
\draw (1,1.8) rectangle (1.6,2.4);
\draw (1.6,0) rectangle (2.2,.6);
\draw (1.6,.6) rectangle (2.2,1.2);
\draw (1.6,1.2) rectangle (2.2,1.8);
\draw (2.2,0) rectangle (2.8,.6);
\draw (2.2,.6) rectangle (2.8,1.2);
\draw (2.8,0) rectangle (3.4,.6);
 \node at (1.3,.3) {5};
 \node at (1.3,.9) {4};
 \node at (1.3,1.5) {2};
 \node at (1.3,2.1) {3};
 \node at (1.9,.3) {7};
 \node at (1.9,.9) {6};
 \node at (1.9,1.5) {1};
 \node at (2.5,.3) {9};
 \node at (2.5,.9) {8};
 \node at (3.1,.3) {10};
 \end{tikzpicture}
\end{center}
\caption{Two standard balanced tableaux.}\label{rothe2}
\end{figure}
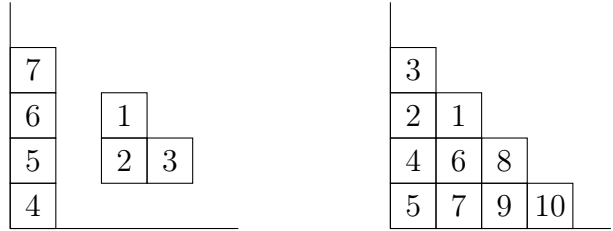

A labelling, with no repeats, of the cells of the Rothe diagram for a permutation  $w\in \ss_n$ with the positive integers in $[\ell(w)]$ is called a {\it standard balanced tableaux} if  for any entry of the diagram, the number of entries to its right in the same row that are greater than it is equal to the number of entries above it in the same column that are smaller.
 Figure \ref{rothe2} shows two standard balanced tableaux on $\mathbb{D}(25431)$ and
 $\mathbb{D}(54321)$.

 Denote the set of all standard balanced tableaux on $\mathbb{D}(w)$ by $\SBT(w)$.

 Let $w\in \ss_n$ be a permutation of length $\ell$ and $a=a_1a_2\cdots a_{\ell}$ a reduced decomposition of $w$.
 Note that each $a_i$ corresponds to a unique inversion in $w$, namely the pair of numbers transposed by $a_i$ in the product
 $$w=s_{a_1}s_{a_2}\cdots s_{a_{\ell}},$$
 and $a$ is determined uniquely by the order in which these inversions are carried out. Since $(p,q)\in \mathbb{D}(w)$ if and only if $(p,q)$ is an inversion of $w$, we
 define the labelling  $P_a$ of $\mathbb{D}(w)$ by setting $P_a(p,q)=i$ if $a_i$ transposes $p$ and $q$,
 where $p>q$.
 The tableau
 $P_a$ is a standard balanced tableau called the {\it canonical labelling} of $\mathbb{D}(w)$ induced by $a$.

 Fomin {\it et al.} \cite{reiner} proved that the map $a\mapsto P_a$ defines a bijection between $\Red(w)$ and $\SBT(w)$.

\begin{example}\label{ex:action}
 Let $w=4321\in S_4$ and let $a=213213\in \Red(w)$ be a reduced decomposition of $w$. The action of each $a_i$ is illustrated in the following table, from which we get the standard balanced tableau $D_a$.

\vspace*{5mm}

\parbox{8cm}{
 %\begin{table}[htp]
\begin{center}
\begin{tabular}{c|cccc}
&1&2&3&4\\
\hline
$a_1=2$&1&{\bf 3}&{\bf 2}&4\\
$a_2=1$&{\bf 3}&{\bf 1}&2&4\\
$a_3=3$&3&1&{\bf 4}&{\bf 2}\\
$a_4=2$&3&{\bf 4}&{\bf 1}&2\\
$a_5=1$&{\bf 4}&{\bf 3}&1&2\\
$a_6=3$&4&3&{\bf 2}&{\bf 1}
\end{tabular}
\end{center}
%\end{table}
}
\parbox{8cm}{
 \begin{tikzpicture}
 \node at (0,1) {$P_a$=};
 \draw (1,0) rectangle (1.6,.6);
\draw (1,.6) rectangle (1.6,1.2);
\draw (1,1.2) rectangle (1.6,1.8);
\draw (1.6,0) rectangle (2.2,.6);
\draw (1.6,.6) rectangle (2.2,1.2);
\draw (2.2,0) rectangle (2.8,.6);
 \node at (1.3,.3) {4};
 \node at (1.3,.9) {2};
 \node at (1.3,1.5) {6};
 \node at (1.9,.3) {3};
 \node at (1.9,.9) {1};
 \node at (2.5,.3) {5};
 \node at (1.3,-.3) {\tiny 1};
 \node at (1.9,-.3) {\tiny 2};
 \node at (2.5,-.3) {\tiny 3};
 \node at (3.1,-.3) {\tiny 4};
 \node at (.7,.3) {\tiny 4};
 \node at (.7,.9) {\tiny 3};
 \node at (.7,1.5) {\tiny 2};
 \node at (.7,2.1) {\tiny 1};
 \end{tikzpicture}
 }
\end{example}

 \section{The commutation graph of a permutation}

 We define a relation $\sim$ on the set $\Red(w)$  by setting $a\sim b$ if and only if $a$ and $b$ differ by a sequence of commutations. This is an equivalence relation and the classes it defines are the commutation classes of $w$, denoted by $[a]$. We write $a\sbr b$ (resp. $a\lbr b$) when $a$ and $b$ differ by a single commutation  (resp. long braid relation), and $[a]\lbr[b]$ when those classes differ by a single long braid relation, {\it i.e.} when exist $a'\in[a]$ and $b'\in[b]$ such that $a'\lbr b'$.

\begin{definition}
 Let $w\in\ss_n$. The {\it commutation graph} $C(w)$ has vertex set the commutation classes of $\Red(w)$, and an edge connecting two classes when they differ by a long braid relation.
\end{definition}

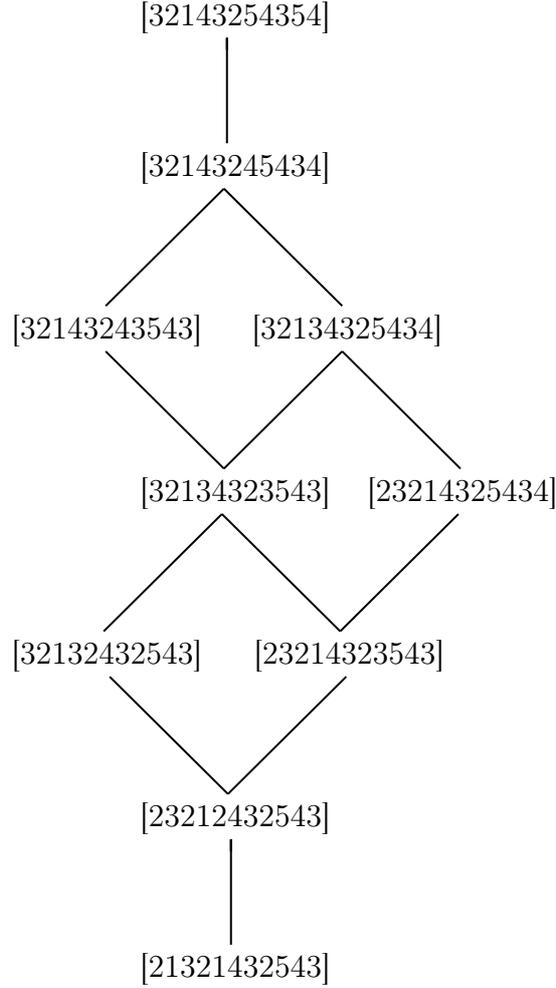
\begin{figure}
  \centering

\tikzset{every picture/.style={line width=0.75pt}} %set default line width to 0.75pt

\begin{tikzpicture}[x=0.75pt,y=0.75pt,yscale=-1,xscale=1]
%uncomment if require: \path (0,1304); %set diagram left start at 0, and has height of 1304

%Straight Lines [id:da052757371332554204]
\draw    (312.5,511.12) -- (312.5,505.12) -- (312.5,558.12) ;
%Straight Lines [id:da28027425234967573]
\draw    (252,423.28) -- (311,482.28) ;
%Straight Lines [id:da12660675503722496]
\draw    (370,423.28) -- (311,482.28) ;
%Straight Lines [id:da22138593471172197]
\draw    (308,341.44) -- (249,400.44) ;
%Straight Lines [id:da4156798435906517]
\draw    (308,341.44) -- (367,400.44) ;
%Straight Lines [id:da2500237313886833]
\draw    (426,341.44) -- (367,400.44) ;
%Straight Lines [id:da7821236222193224]
\draw    (250,259.6) -- (309,318.6) ;
%Straight Lines [id:da6663325689680735]
\draw    (368,259.6) -- (309,318.6) ;
%Straight Lines [id:da031215834112047647]
\draw    (368,259.6) -- (427,318.6) ;
%Straight Lines [id:da6495204833515822]
\draw    (309,177.76) -- (250,236.76) ;
%Straight Lines [id:da7980184777210217]
\draw    (309,177.76) -- (368,236.76) ;
%Straight Lines [id:da32141917832697886]
\draw    (310.5,107.92) -- (310.5,101.92) -- (310.5,154.92) ;
%Straight Lines [id:da12029316447279959]
%\draw    (303,671.5) -- (381,671.5) ;
%Straight Lines [id:da8789487589863081]
%\draw    (412,708.5) -- (412,747.5) ;

\draw (265.5,561) node [anchor=north west][inner sep=0.75pt]   [align=left] {[21321432543]};
% Text Node
\draw (265.5,485.2) node [anchor=north west][inner sep=0.75pt]   [align=left] {[23212432543]};
% Text Node
\draw (265.5,157.84) node [anchor=north west][inner sep=0.75pt]   [align=left] {[32143245434]};
% Text Node
\draw (265.5,82) node [anchor=north west][inner sep=0.75pt]   [align=left] {[32143254354]};
% Text Node
\draw (201.5,239.68) node [anchor=north west][inner sep=0.75pt]   [align=left] {[32143243543]};
% Text Node
\draw (321.5,239.68) node [anchor=north west][inner sep=0.75pt]   [align=left] {[32134325434]};
% Text Node
\draw (265.5,321.52) node [anchor=north west][inner sep=0.75pt]   [align=left] {[32134323543]};
% Text Node
\draw (379,321.52) node [anchor=north west][inner sep=0.75pt]   [align=left] {[23214325434]};
% Text Node
\draw (201.5,403.36) node [anchor=north west][inner sep=0.75pt]   [align=left] {[32132432543]};
% Text Node
\draw (322.5,403.36) node [anchor=north west][inner sep=0.75pt]   [align=left] {[23214323543]};

\end{tikzpicture}

 \caption{The graph $C(456312)$.}\label{fig:grafo}
\end{figure}

Note that since $C(w)$ can be obtained from $G(w)$ by contracting commutation edges, the connectivity of $G(w)$ implies the connectivity of $C(w)$.
The {\it distance} $d([a],[b])$ between two commutation classes $[a]$ and $[b]$ on $C(w)$ is the length of a shortest path joining $[a]$ and $[b]$, that is the number of edges in such path.
The {\it eccentricity} of a class $[a]$ is the distance to a farthest
commutative class from $[a]$. The {\it radius} and {\it diameter} of $C(w)$ are the minimum and
maximum eccentricities, respectively.
Figure \ref{fig:grafo} shows the commutation graph $C(456312)$, which has diameter 6 and radius 3.

Next, following \cite{sami}, we define analogs of commutation and braid relations for balanced tableaux.

\begin{definition}
 Given a permutation $w$ and an integer $1\leq i< \ell(w)$, the map $\cc_i$ acts on the tableaux in $\SBT(w)$
 for which the labels $i$ and $i+1$ are not in the same row nor the same column,
 by interchanging the labels $i$ and $i+1$.
\end{definition}

 It easy to check that the map $\cc_i$ is  well-defined, since for $R\in\SBT(w)$, if $i$ and $i+1$ are not in the same row nor the same column, then interchanging them keep the balanced condition true since all other entries compare the same with $i$ and with $i+1$. Thus, $\cc_i(R)\in\SBT(w)$. The map $\cc_i$ is clearly an involution.

 Let $a=pa_ia_{i+1}q$ be a reduced decompositions of $w\in \ss_n$ with $|a_i-a_{i+1}|>1$. Then, $a_i$ and $a_{i+1}$ correspond to inversions on disjoint set of integers, say $(y,x)$ and $(k,z)$, with $x<y$ and $z<k$. This means that the cells $(y,x)$ and $(k,w)$ of tableau $P_a$ are in distinct rows and columns. Swapping the labels  $i$ and $i+1$ of these cells gives the tableau $P_{a'}$, where $a'=pa'_ia'_{i+1}q\sbr a$, with $a'_i=a_{i+1}$ and $a'_{i+1}=a_i$. That is, we have $a\sbr a'$ if and only if $\cc_i(P_a)=P_{a'}$.

 \begin{definition}
 Given a permutation $w$ and an integer $1<i< \ell(w)$, the map $\bb_i$ acts on the tableaux in $\SBT(w)$
 having one of the  labels $i-1$ or $i+1$ in the same column and above $i$ and the other  in the same row and right of $i$,
  by interchanging  the labels $i-1$ and $i+1$.
\end{definition}

 For $R\in \SBT(w)$, if $i\pm 1$ is in the same row as $i$, and $i\mp 1$ is in the same column, then swapping them maintains the tableaux balanced since  all entries  $j\notin \{i-1,i,i+i\}$   compares the same with both integers $i\pm1$. Thus, the map $\bb_i$ is an  involution.

 Let $a=pa_{i-1}a_ia_{i+1}q$ be a reduced decomposition of $w\in \ss_n$, with $a_{i-1}=a_{i+1}=a_i\pm 1$.
 Then, $a_{i-1}, a_i$ and $a_{i+1}$ correspond to inversions  $(y,x)$, $(z,x)$ and $(z,y)$, with $x<y<z$. %Note that in the tableau $P_a$, the length of the hook with end cells $(z,y)$ and $(y,x)$ is $z-x+1+|w_z-w_y|$.
 Interchanging the labels $i+1$ and $i-1$ in $P_a$ gives the tableau $P_{a'}$, where $a=pa'_{i-1}a'_ia'_{i+1}q$
 with $a'_{i-1}=a'_{i+1}=a_i$ and $a'_i=a_{i+1}$, showing that $a\lbr a'$ if and only if $\bb_i(P_a)= P_{a'}$.

 Figure \ref{fig:grafo1} shows the action of the maps $\cc_i$ and $\bb_j$ on some tableaux in $\SBT(456312)$ for the classes of the bottom three levels of graph $C(456312)$ depicted in Figure \ref{fig:grafo}. The tableaux in bold font correspond to the words represented in the graph, and the dashed lines correspond to maps that can also act in the tableaux.

 \begin{figure}

 \tikzset{every picture/.style={line width=0.75pt}} %set default line width to 0.75pt

\tikzset{every picture/.style={line width=0.75pt}} %set default line width to 0.75pt

\begin{tikzpicture}[x=0.75pt,y=0.75pt,yscale=-1,xscale=1]
%uncomment if require: \path (0,1304); %set diagram left start at 0, and has height of 1304
\footnotesize{
%Shape: Rectangle [id:dp15826189645137045]
\draw   (261.26,1132.1) -- (277.41,1132.1) -- (277.41,1148.25) -- (261.26,1148.25) -- cycle ;
%Shape: Rectangle [id:dp2884390225986879]
\draw   (261.26,1148.25) -- (277.41,1148.25) -- (277.41,1164.4) -- (261.26,1164.4) -- cycle ;
%Shape: Square [id:dp6298547687090443]
\draw   (261.26,1164.4) -- (277.41,1164.4) -- (277.41,1180.55) -- (261.26,1180.55) -- cycle ;
%Shape: Rectangle [id:dp11617397765245108]
\draw   (261.26,1180.55) -- (277.41,1180.55) -- (277.41,1196.7) -- (261.26,1196.7) -- cycle ;
%Shape: Square [id:dp984311641503655]
\draw   (277.41,1164.4) -- (293.56,1164.4) -- (293.56,1180.55) -- (277.41,1180.55) -- cycle ;
%Shape: Rectangle [id:dp09691463574802217]
\draw   (277.41,1148.25) -- (293.56,1148.25) -- (293.56,1164.4) -- (277.41,1164.4) -- cycle ;
%Shape: Square [id:dp9772167322774792]
\draw   (277.41,1180.55) -- (293.56,1180.55) -- (293.56,1196.7) -- (277.41,1196.7) -- cycle ;
%Shape: Rectangle [id:dp6337284115017077]
\draw   (277.41,1132.1) -- (293.56,1132.1) -- (293.56,1148.25) -- (277.41,1148.25) -- cycle ;
%Shape: Square [id:dp5769833997956095]
\draw   (293.56,1164.4) -- (309.71,1164.4) -- (309.71,1180.55) -- (293.56,1180.55) -- cycle ;
%Shape: Rectangle [id:dp041499516264007275]
\draw   (293.56,1148.25) -- (309.71,1148.25) -- (309.71,1164.4) -- (293.56,1164.4) -- cycle ;
%Shape: Square [id:dp4379418982246095]
\draw   (293.56,1180.55) -- (309.71,1180.55) -- (309.71,1196.7) -- (293.56,1196.7) -- cycle ;
%Shape: Rectangle [id:dp4998918935205967]
\draw   (411.26,1132.1) -- (427.41,1132.1) -- (427.41,1148.25) -- (411.26,1148.25) -- cycle ;
%Shape: Rectangle [id:dp012423776799128872]
\draw   (411.26,1148.25) -- (427.41,1148.25) -- (427.41,1164.4) -- (411.26,1164.4) -- cycle ;
%Shape: Square [id:dp43404174772828985]
\draw   (411.26,1164.4) -- (427.41,1164.4) -- (427.41,1180.55) -- (411.26,1180.55) -- cycle ;
%Shape: Rectangle [id:dp15889147139070858]
\draw   (411.26,1180.55) -- (427.41,1180.55) -- (427.41,1196.7) -- (411.26,1196.7) -- cycle ;
%Shape: Square [id:dp9539547802524033]
\draw   (427.41,1164.4) -- (443.56,1164.4) -- (443.56,1180.55) -- (427.41,1180.55) -- cycle ;
%Shape: Rectangle [id:dp7229304786504493]
\draw   (427.41,1148.25) -- (443.56,1148.25) -- (443.56,1164.4) -- (427.41,1164.4) -- cycle ;
%Shape: Square [id:dp8676365107436523]
\draw   (427.41,1180.55) -- (443.56,1180.55) -- (443.56,1196.7) -- (427.41,1196.7) -- cycle ;
%Shape: Rectangle [id:dp7367309162341649]
\draw   (427.41,1132.1) -- (443.56,1132.1) -- (443.56,1148.25) -- (427.41,1148.25) -- cycle ;
%Shape: Square [id:dp5839135609501374]
\draw   (443.56,1164.4) -- (459.71,1164.4) -- (459.71,1180.55) -- (443.56,1180.55) -- cycle ;
%Shape: Rectangle [id:dp556793410060278]
\draw   (443.56,1148.25) -- (459.71,1148.25) -- (459.71,1164.4) -- (443.56,1164.4) -- cycle ;
%Shape: Square [id:dp24296571582514437]
\draw   (443.56,1180.55) -- (459.71,1180.55) -- (459.71,1196.7) -- (443.56,1196.7) -- cycle ;
%Straight Lines [id:da12029316447279959]
\draw    (323,1164.5) -- (401,1164.5) ;
%Shape: Rectangle [id:dp1054381737444583]
\draw   (410.26,1005.1) -- (426.41,1005.1) -- (426.41,1021.25) -- (410.26,1021.25) -- cycle ;
%Shape: Rectangle [id:dp9254975678966999]
\draw   (410.26,1021.25) -- (426.41,1021.25) -- (426.41,1037.4) -- (410.26,1037.4) -- cycle ;
%Shape: Square [id:dp576904225045058]
\draw   (410.26,1037.4) -- (426.41,1037.4) -- (426.41,1053.55) -- (410.26,1053.55) -- cycle ;
%Shape: Rectangle [id:dp8621575403757515]
\draw   (410.26,1053.55) -- (426.41,1053.55) -- (426.41,1069.7) -- (410.26,1069.7) -- cycle ;
%Shape: Square [id:dp46391895127253413]
\draw   (426.41,1037.4) -- (442.56,1037.4) -- (442.56,1053.55) -- (426.41,1053.55) -- cycle ;
%Shape: Rectangle [id:dp892526066880212]
\draw   (426.41,1021.25) -- (442.56,1021.25) -- (442.56,1037.4) -- (426.41,1037.4) -- cycle ;
%Shape: Square [id:dp46319686070576616]
\draw   (426.41,1053.55) -- (442.56,1053.55) -- (442.56,1069.7) -- (426.41,1069.7) -- cycle ;
%Shape: Rectangle [id:dp2783940833398586]
\draw   (426.41,1005.1) -- (442.56,1005.1) -- (442.56,1021.25) -- (426.41,1021.25) -- cycle ;
%Shape: Square [id:dp6085291700879749]
\draw   (442.56,1037.4) -- (458.71,1037.4) -- (458.71,1053.55) -- (442.56,1053.55) -- cycle ;
%Shape: Rectangle [id:dp3588719998135983]
\draw   (442.56,1021.25) -- (458.71,1021.25) -- (458.71,1037.4) -- (442.56,1037.4) -- cycle ;
%Shape: Square [id:dp6845726293247381]
\draw   (442.56,1053.55) -- (458.71,1053.55) -- (458.71,1069.7) -- (442.56,1069.7) -- cycle ;
%Straight Lines [id:da8789487589863081]
\draw    (429,1080.5) -- (429,1119.5) ;
%Straight Lines [id:da25872517303129716]
\draw[dashed]    (470,1150) -- (520,1150) ;
%Straight Lines [id:da3070521679086349]
\draw[dashed]    (470,1180) -- (520,1180) ;
%Straight Lines [id:da6193660275050976]
\draw[dashed]    (200,1150) -- (250,1150) ;
%Straight Lines [id:da9315087361035548]
\draw[dashed]    (200,1180) -- (250,1180) ;
%Straight Lines [id:da7838185262939261]
\draw    (322,1039.5) -- (400,1039.5) ;
%Shape: Rectangle [id:dp04906053103908614]
\draw   (263.26,1006.1) -- (279.41,1006.1) -- (279.41,1022.25) -- (263.26,1022.25) -- cycle ;
%Shape: Rectangle [id:dp7046235121264421]
\draw   (263.26,1022.25) -- (279.41,1022.25) -- (279.41,1038.4) -- (263.26,1038.4) -- cycle ;
%Shape: Square [id:dp7205489876663878]
\draw   (263.26,1038.4) -- (279.41,1038.4) -- (279.41,1054.55) -- (263.26,1054.55) -- cycle ;
%Shape: Rectangle [id:dp45944708284791624]
\draw   (263.26,1054.55) -- (279.41,1054.55) -- (279.41,1070.7) -- (263.26,1070.7) -- cycle ;
%Shape: Square [id:dp42411692509130705]
\draw   (279.41,1038.4) -- (295.56,1038.4) -- (295.56,1054.55) -- (279.41,1054.55) -- cycle ;
%Shape: Rectangle [id:dp38052540663681955]
\draw   (279.41,1022.25) -- (295.56,1022.25) -- (295.56,1038.4) -- (279.41,1038.4) -- cycle ;
%Shape: Square [id:dp9071875688487594]
\draw   (279.41,1054.55) -- (295.56,1054.55) -- (295.56,1070.7) -- (279.41,1070.7) -- cycle ;
%Shape: Rectangle [id:dp7337649103995842]
\draw   (279.41,1006.1) -- (295.56,1006.1) -- (295.56,1022.25) -- (279.41,1022.25) -- cycle ;
%Shape: Square [id:dp860666937984145]
\draw   (295.56,1038.4) -- (311.71,1038.4) -- (311.71,1054.55) -- (295.56,1054.55) -- cycle ;
%Shape: Rectangle [id:dp6485014863004397]
\draw   (295.56,1022.25) -- (311.71,1022.25) -- (311.71,1038.4) -- (295.56,1038.4) -- cycle ;
%Shape: Square [id:dp16583041108148122]
\draw   (295.56,1054.55) -- (311.71,1054.55) -- (311.71,1070.7) -- (295.56,1070.7) -- cycle ;
%Straight Lines [id:da9187945719299311]
\draw[dashed]    (200,1023) -- (250,1023) ;
%Straight Lines [id:da3318249550695491]
\draw[dashed]    (200,1053) -- (250,1053) ;
%Straight Lines [id:da9609895776757409]
\draw[dashed]    (470,1039.5) -- (520,1039.5) ;
%Shape: Rectangle [id:dp5092874838104042]
\draw   (412.26,876.1) -- (428.41,876.1) -- (428.41,892.25) -- (412.26,892.25) -- cycle ;
%Shape: Rectangle [id:dp3196211410985135]
\draw   (412.26,892.25) -- (428.41,892.25) -- (428.41,908.4) -- (412.26,908.4) -- cycle ;
%Shape: Square [id:dp3467702322952573]
\draw   (412.26,908.4) -- (428.41,908.4) -- (428.41,924.55) -- (412.26,924.55) -- cycle ;
%Shape: Rectangle [id:dp9360284686714853]
\draw   (412.26,924.55) -- (428.41,924.55) -- (428.41,940.7) -- (412.26,940.7) -- cycle ;
%Shape: Square [id:dp050519076970731325]
\draw   (428.41,908.4) -- (444.56,908.4) -- (444.56,924.55) -- (428.41,924.55) -- cycle ;
%Shape: Rectangle [id:dp025196778126689745]
\draw   (428.41,892.25) -- (444.56,892.25) -- (444.56,908.4) -- (428.41,908.4) -- cycle ;
%Shape: Square [id:dp43839674939109585]
\draw   (428.41,924.55) -- (444.56,924.55) -- (444.56,940.7) -- (428.41,940.7) -- cycle ;
%Shape: Rectangle [id:dp18499104355969598]
\draw   (428.41,876.1) -- (444.56,876.1) -- (444.56,892.25) -- (428.41,892.25) -- cycle ;
%Shape: Square [id:dp8251081166585044]
\draw   (444.56,908.4) -- (460.71,908.4) -- (460.71,924.55) -- (444.56,924.55) -- cycle ;
%Shape: Rectangle [id:dp26524075843818773]
\draw   (444.56,892.25) -- (460.71,892.25) -- (460.71,908.4) -- (444.56,908.4) -- cycle ;
%Shape: Square [id:dp7032094173361565]
\draw   (444.56,924.55) -- (460.71,924.55) -- (460.71,940.7) -- (444.56,940.7) -- cycle ;
%Straight Lines [id:da30325573505645487]
\draw    (471,907.5) -- (549,907.5) ;
%Shape: Rectangle [id:dp6488555731120487]
\draw   (559.26,875.1) -- (575.41,875.1) -- (575.41,891.25) -- (559.26,891.25) -- cycle ;
%Shape: Rectangle [id:dp780209407763605]
\draw   (559.26,891.25) -- (575.41,891.25) -- (575.41,907.4) -- (559.26,907.4) -- cycle ;
%Shape: Square [id:dp6129180180299354]
\draw   (559.26,907.4) -- (575.41,907.4) -- (575.41,923.55) -- (559.26,923.55) -- cycle ;
%Shape: Rectangle [id:dp8984789537815259]
\draw   (559.26,923.55) -- (575.41,923.55) -- (575.41,939.7) -- (559.26,939.7) -- cycle ;
%Shape: Square [id:dp44213593554270747]
\draw   (575.41,907.4) -- (591.56,907.4) -- (591.56,923.55) -- (575.41,923.55) -- cycle ;
%Shape: Rectangle [id:dp6906491703763085]
\draw   (575.41,891.25) -- (591.56,891.25) -- (591.56,907.4) -- (575.41,907.4) -- cycle ;
%Shape: Square [id:dp5471091633344833]
\draw   (575.41,923.55) -- (591.56,923.55) -- (591.56,939.7) -- (575.41,939.7) -- cycle ;
%Shape: Rectangle [id:dp32921116450968846]
\draw   (575.41,875.1) -- (591.56,875.1) -- (591.56,891.25) -- (575.41,891.25) -- cycle ;
%Shape: Square [id:dp36679286029216085]
\draw   (591.56,907.4) -- (607.71,907.4) -- (607.71,923.55) -- (591.56,923.55) -- cycle ;
%Shape: Rectangle [id:dp5992084223629506]
\draw   (591.56,891.25) -- (607.71,891.25) -- (607.71,907.4) -- (591.56,907.4) -- cycle ;
%Shape: Square [id:dp4705022261208369]
\draw   (591.56,923.55) -- (607.71,923.55) -- (607.71,939.7) -- (591.56,939.7) -- cycle ;
%Straight Lines [id:da7536429239397286]
\draw[dashed]    (355,893.5) -- (403,893.5) ;
%Straight Lines [id:da4537224898306962]
\draw[dashed]    (354,919.5) -- (402,919.5) ;
%Straight Lines [id:da5628908239594466]
\draw[dashed]    (620,895.5) -- (668,895.5) ;
%Straight Lines [id:da9248664883483442]
\draw[dashed]    (619,921.5) -- (667,921.5) ;
%Straight Lines [id:da14487893797243911]
\draw    (431,950.5) -- (431,989.5) ;
%Straight Lines [id:da4450870203633579]
\draw    (241,950.5) -- (279,991.5) ;
%Shape: Rectangle [id:dp6399490453111154]
\draw   (181.26,876.1) -- (197.41,876.1) -- (197.41,892.25) -- (181.26,892.25) -- cycle ;
%Shape: Rectangle [id:dp48349449863112803]
\draw   (181.26,892.25) -- (197.41,892.25) -- (197.41,908.4) -- (181.26,908.4) -- cycle ;
%Shape: Square [id:dp667589940805968]
\draw   (181.26,908.4) -- (197.41,908.4) -- (197.41,924.55) -- (181.26,924.55) -- cycle ;
%Shape: Rectangle [id:dp5429508796971256]
\draw   (181.26,924.55) -- (197.41,924.55) -- (197.41,940.7) -- (181.26,940.7) -- cycle ;
%Shape: Square [id:dp9906765061988458]
\draw   (197.41,908.4) -- (213.56,908.4) -- (213.56,924.55) -- (197.41,924.55) -- cycle ;
%Shape: Rectangle [id:dp623280037626172]
\draw   (197.41,892.25) -- (213.56,892.25) -- (213.56,908.4) -- (197.41,908.4) -- cycle ;
%Shape: Square [id:dp29064502215664634]
\draw   (197.41,924.55) -- (213.56,924.55) -- (213.56,940.7) -- (197.41,940.7) -- cycle ;
%Shape: Rectangle [id:dp5635621872152676]
\draw   (197.41,876.1) -- (213.56,876.1) -- (213.56,892.25) -- (197.41,892.25) -- cycle ;
%Shape: Square [id:dp2875090333263215]
\draw   (213.56,908.4) -- (229.71,908.4) -- (229.71,924.55) -- (213.56,924.55) -- cycle ;
%Shape: Rectangle [id:dp7666036028529348]
\draw   (213.56,892.25) -- (229.71,892.25) -- (229.71,908.4) -- (213.56,908.4) -- cycle ;
%Shape: Square [id:dp5677902831246369]
\draw   (213.56,924.55) -- (229.71,924.55) -- (229.71,940.7) -- (213.56,940.7) -- cycle ;
%Straight Lines [id:da5117661009471643]
\draw[dashed]    (125,894.5) -- (173,894.5) ;
%Straight Lines [id:da3954273717822001]
\draw[dashed]    (124,920.5) -- (172,920.5) ;

% Text Node
\draw (263.83,1182.13) node [anchor=north west][inner sep=0.75pt]   [align=left] {\textbf{4}};
% Text Node
\draw (263.83,1133.68) node [anchor=north west][inner sep=0.75pt]   [align=left] {\textbf{2}};
% Text Node
\draw (259.83,1149.83) node [anchor=north west][inner sep=0.75pt]   [align=left] {\textbf{10}};
% Text Node
\draw (263.83,1165.98) node [anchor=north west][inner sep=0.75pt]   [align=left] {\textbf{7}};
% Text Node
\draw (279.98,1149.83) node [anchor=north west][inner sep=0.75pt]   [align=left] {\textbf{9}};
% Text Node
\draw (279.98,1165.98) node [anchor=north west][inner sep=0.75pt]   [align=left] {\textbf{6}};
% Text Node
\draw (279.98,1182.13) node [anchor=north west][inner sep=0.75pt]   [align=left] {\textbf{3}};
% Text Node
\draw (350.85,1150) node [anchor=north west][inner sep=0.75pt]    {$\cc_{2}$};
% Text Node
\draw (279.98,1133.68) node [anchor=north west][inner sep=0.75pt]   [align=left] {\textbf{1}};
% Text Node
\draw (291.13,1149.83) node [anchor=north west][inner sep=0.75pt]   [align=left] {\textbf{11}};
% Text Node
\draw (296.13,1165.98) node [anchor=north west][inner sep=0.75pt]   [align=left] {\textbf{8}};
% Text Node
\draw (296.13,1182.13) node [anchor=north west][inner sep=0.75pt]   [align=left] {\textbf{5}};
% Text Node
\draw (413.83,1182.13) node [anchor=north west][inner sep=0.75pt]   [align=left] {4};
% Text Node
\draw (413.83,1133.68) node [anchor=north west][inner sep=0.75pt]   [align=left] {3};
% Text Node
\draw (409.83,1149.83) node [anchor=north west][inner sep=0.75pt]   [align=left] {10};
% Text Node
\draw (413.83,1165.98) node [anchor=north west][inner sep=0.75pt]   [align=left] {7};
% Text Node
\draw (429.98,1149.83) node [anchor=north west][inner sep=0.75pt]   [align=left] {9};
% Text Node
\draw (429.98,1165.98) node [anchor=north west][inner sep=0.75pt]   [align=left] {6};
% Text Node
\draw (429.98,1182.13) node [anchor=north west][inner sep=0.75pt]   [align=left] {2};
% Text Node
\draw (429.98,1133.68) node [anchor=north west][inner sep=0.75pt]   [align=left] {1};
% Text Node
\draw (441.13,1149.83) node [anchor=north west][inner sep=0.75pt]   [align=left] {11};
% Text Node
\draw (446.13,1165.98) node [anchor=north west][inner sep=0.75pt]   [align=left] {8};
% Text Node
\draw (446.13,1182.13) node [anchor=north west][inner sep=0.75pt]   [align=left] {5};
% Text Node
\draw (412.83,1055.13) node [anchor=north west][inner sep=0.75pt]   [align=left] {\textbf{4}};
% Text Node
\draw (412.83,1006.68) node [anchor=north west][inner sep=0.75pt]   [align=left] {\textbf{5}};
% Text Node
\draw (408.83,1022.83) node [anchor=north west][inner sep=0.75pt]   [align=left] {\textbf{10}};
% Text Node
\draw (412.83,1038.98) node [anchor=north west][inner sep=0.75pt]   [align=left] {\textbf{7}};
% Text Node
\draw (428.98,1022.83) node [anchor=north west][inner sep=0.75pt]   [align=left] {\textbf{9}};
% Text Node
\draw (428.98,1038.98) node [anchor=north west][inner sep=0.75pt]   [align=left] {\textbf{6}};
% Text Node
\draw (428.98,1055.13) node [anchor=north west][inner sep=0.75pt]   [align=left] {\textbf{2}};
% Text Node
\draw (428.98,1006.68) node [anchor=north west][inner sep=0.75pt]   [align=left] {\textbf{1}};
% Text Node
\draw (440.13,1022.83) node [anchor=north west][inner sep=0.75pt]   [align=left] {\textbf{11}};
% Text Node
\draw (445.13,1038.98) node [anchor=north west][inner sep=0.75pt]   [align=left] {\textbf{8}};
% Text Node
\draw (445.13,1055.13) node [anchor=north west][inner sep=0.75pt]   [align=left] {\textbf{3}};
% Text Node
\draw (430,1091.25) node [anchor=north west][inner sep=0.75pt]    {$\bb_{4}$};
% Text Node
\draw (485,1135) node [anchor=north west][inner sep=0.75pt]    {$\cc_{5}$};
% Text Node
\draw (485,1180) node [anchor=north west][inner sep=0.75pt]    {$\cc_{8}$};
% Text Node
\draw (215,1135) node [anchor=north west][inner sep=0.75pt]    {$\cc_{5}$};
% Text Node
\draw (215,1180) node [anchor=north west][inner sep=0.75pt]    {$\cc_{8}$};
% Text Node
\draw (350,1025) node [anchor=north west][inner sep=0.75pt]    {$\cc_{5}$};
% Text Node
\draw (265.83,1056.13) node [anchor=north west][inner sep=0.75pt]   [align=left] {4};
% Text Node
\draw (265.83,1007.68) node [anchor=north west][inner sep=0.75pt]   [align=left] {6};
% Text Node
\draw (261.83,1023.83) node [anchor=north west][inner sep=0.75pt]   [align=left] {10};
% Text Node
\draw (265.83,1039.98) node [anchor=north west][inner sep=0.75pt]   [align=left] {7};
% Text Node
\draw (281.98,1023.83) node [anchor=north west][inner sep=0.75pt]   [align=left] {9};
% Text Node
\draw (281.98,1039.98) node [anchor=north west][inner sep=0.75pt]   [align=left] {5};
% Text Node
\draw (281.98,1056.13) node [anchor=north west][inner sep=0.75pt]   [align=left] {2};
% Text Node
\draw (281.98,1007.68) node [anchor=north west][inner sep=0.75pt]   [align=left] {1};
% Text Node
\draw (293.13,1023.83) node [anchor=north west][inner sep=0.75pt]   [align=left] {11};
% Text Node
\draw (298.13,1039.98) node [anchor=north west][inner sep=0.75pt]   [align=left] {8};
% Text Node
\draw (298.13,1056.13) node [anchor=north west][inner sep=0.75pt]   [align=left] {3};
% Text Node
\draw (215,1010) node [anchor=north west][inner sep=0.75pt]    {$\cc_{4}$};
% Text Node
\draw (215,1055) node [anchor=north west][inner sep=0.75pt]    {$\cc_{8}$};
% Text Node
\draw (485,1025) node [anchor=north west][inner sep=0.75pt]    {$\cc_{8}$};
% Text Node
\draw (414.83,926.13) node [anchor=north west][inner sep=0.75pt]   [align=left] {4};
% Text Node
\draw (414.83,877.68) node [anchor=north west][inner sep=0.75pt]   [align=left] {5};
% Text Node
\draw (410.83,893.83) node [anchor=north west][inner sep=0.75pt]   [align=left] {10};
% Text Node
\draw (414.83,909.98) node [anchor=north west][inner sep=0.75pt]   [align=left] {7};
% Text Node
\draw (430.98,893.83) node [anchor=north west][inner sep=0.75pt]   [align=left] {9};
% Text Node
\draw (430.98,909.98) node [anchor=north west][inner sep=0.75pt]   [align=left] {6};
% Text Node
\draw (430.98,926.13) node [anchor=north west][inner sep=0.75pt]   [align=left] {2};
% Text Node
\draw (430.98,877.68) node [anchor=north west][inner sep=0.75pt]   [align=left] {3};
% Text Node
\draw (442.13,893.83) node [anchor=north west][inner sep=0.75pt]   [align=left] {11};
% Text Node
\draw (447.13,909.98) node [anchor=north west][inner sep=0.75pt]   [align=left] {8};
% Text Node
\draw (447.13,926.13) node [anchor=north west][inner sep=0.75pt]   [align=left] {1};
% Text Node
\draw (498.85,894.25) node [anchor=north west][inner sep=0.75pt]    {$\cc_{3}$};
% Text Node
\draw (561.83,925.13) node [anchor=north west][inner sep=0.75pt]   [align=left] {\textbf{3}};
% Text Node
\draw (561.83,876.68) node [anchor=north west][inner sep=0.75pt]   [align=left] {\textbf{5}};
% Text Node
\draw (557.83,892.83) node [anchor=north west][inner sep=0.75pt]   [align=left] {\textbf{10}};
% Text Node
\draw (561.83,908.98) node [anchor=north west][inner sep=0.75pt]   [align=left] {\textbf{7}};
% Text Node
\draw (577.98,892.83) node [anchor=north west][inner sep=0.75pt]   [align=left] {\textbf{9}};
% Text Node
\draw (577.98,908.98) node [anchor=north west][inner sep=0.75pt]   [align=left] {\textbf{6}};
% Text Node
\draw (577.98,925.13) node [anchor=north west][inner sep=0.75pt]   [align=left] {\textbf{2}};
% Text Node
\draw (577.98,876.68) node [anchor=north west][inner sep=0.75pt]   [align=left] {\textbf{4}};
% Text Node
\draw (589.13,892.83) node [anchor=north west][inner sep=0.75pt]   [align=left] {\textbf{11}};
% Text Node
\draw (594.13,908.98) node [anchor=north west][inner sep=0.75pt]   [align=left] {\textbf{8}};
% Text Node
\draw (594.13,925.13) node [anchor=north west][inner sep=0.75pt]   [align=left] {\textbf{1}};
% Text Node
\draw (367.85,880.25) node [anchor=north west][inner sep=0.75pt]    {$\cc_{5}$};
% Text Node
\draw (367.5,921.25) node [anchor=north west][inner sep=0.75pt]    {$\cc_{8}$};
% Text Node
\draw (632.85,882.25) node [anchor=north west][inner sep=0.75pt]    {$\cc_{5}$};
% Text Node
\draw (632.5,923.25) node [anchor=north west][inner sep=0.75pt]    {$\cc_{8}$};
% Text Node
\draw (443.85,961.25) node [anchor=north west][inner sep=0.75pt]    {$\bb_{2}$};
% Text Node
\draw (266.85,963.25) node [anchor=north west][inner sep=0.75pt]    {$\bb_{7}$};
% Text Node
\draw (183.83,926.13) node [anchor=north west][inner sep=0.75pt]   [align=left] {\textbf{4}};
% Text Node
\draw (183.83,877.68) node [anchor=north west][inner sep=0.75pt]   [align=left] {\textbf{8}};
% Text Node
\draw (179.83,893.83) node [anchor=north west][inner sep=0.75pt]   [align=left] {\textbf{10}};
% Text Node
\draw (183.83,909.98) node [anchor=north west][inner sep=0.75pt]   [align=left] {\textbf{7}};
% Text Node
\draw (199.98,893.83) node [anchor=north west][inner sep=0.75pt]   [align=left] {\textbf{9}};
% Text Node
\draw (199.98,909.98) node [anchor=north west][inner sep=0.75pt]   [align=left] {\textbf{5}};
% Text Node
\draw (199.98,926.13) node [anchor=north west][inner sep=0.75pt]   [align=left] {\textbf{2}};
% Text Node
\draw (199.98,877.68) node [anchor=north west][inner sep=0.75pt]   [align=left] {\textbf{1}};
% Text Node
\draw (211.13,893.83) node [anchor=north west][inner sep=0.75pt]   [align=left] {\textbf{11}};
% Text Node
\draw (216.13,909.98) node [anchor=north west][inner sep=0.75pt]   [align=left] {\textbf{6}};
% Text Node
\draw (216.13,926.13) node [anchor=north west][inner sep=0.75pt]   [align=left] {\textbf{3}};
% Text Node
\draw (137.85,881.25) node [anchor=north west][inner sep=0.75pt]    {$\cc_{4}$};
% Text Node
\draw (137.5,922.25) node [anchor=north west][inner sep=0.75pt]    {$\cc_{8}$};

%Straight Lines [id:zeluis]
\draw[dashed]    (202,831.5) -- (202,870.5) ;
% Text Node
\draw (204.85,842.25) node [anchor=north west][inner sep=0.75pt]    {$\bb_{2}$};
}
\end{tikzpicture}
  \caption{Examples of the action of maps $\cc_i$ and $\bb_j$.}\label{fig:grafo1}
\end{figure}

 We end this section with a result that will be useful in the sequel. Recall that a consecutive substring of a word is called a factor.

\begin{lemma}\label{lemma:admissiblehook}
Suppose $R\in \SBT(w)$ has labels $i-1$ and $i$ in positions $(z,y)$ and $(z,x)$, respectively,  with $x<y<z$. Then, $R$ cannot have a  label $i+1$ in position $(y',x)$ with $y\neq y'$.
\end{lemma}
\begin{proof}
Let $R=P_a$, with $a=a_1\cdots a_{\ell}$ a reduced decomposition of $w$.
Then, the factor $a_{i-1}a_i$ of $a$ acts on the factor $xyz$ of $a_1\cdots a_{i-2}$ transforming it  into $zxy$. It follows that there cannot be a  label $i+1$ in position $(y'x)$ of $P_a$ with $y\neq y'$, since it does not correspond to an inversion  of integers in consecutive positions of the word $a_1\cdots a_i$.
 \end{proof}

 \section{A statistic on $C(w)$}

 In this section we use the statistic $\tt_w$ on $C(w)$, used extensively on $321$-pattern avoidance problems (see  {\it e.g.} \cite{tenner1,tenner2} and the references therein),  which allows the computation of the diameter of the commutation graph of any permutation of the symmetric group.

 \begin{definition}
 Given a permutation $w\in\ss_n$, define the set $\tt_w$ as the collection of all 321-occurrences in $w$ (see also \cite{Dehornoy,ziegler}), that is the set
 $$\tt_w=\{(w_k,w_j,w_i):\,w_i>w_j>w_k\text{ and }i<j<k\},$$
form by all triples $(x,y,z)$ in $[n]$ such that $(z,y)$, $(z,x)$ and $(y,x)$ are inversions for $w$.
\end{definition}

Note that each triple $(x,y,z)$ in $T_w$ corresponds to the endpoints cells $(z,y)$, $(z,x)$ and $(y,x)$ of a hook in the Rothe diagram of $w$.
 For instance, for the permutation $456312\in S_6$ we have
$$\tt_{456312}=\{(1,3,4),(1,3,5),(1,3,6),(2,3,4),(2,3,5),(2,3,6)\},$$
while the set $\tt_{\wo}$, for the longest permutation $\wo=n\cdots 21$ of $\ss_n$  is formed by all $\binom{n}{3}$ triples $(x,y,z)$.

\begin{definition}
 Given a permutation $w\in\ss_n$, define the map $\Gamma$ on the cartesian product  $\Red(w)\times T_w$ by setting
 $$\Gamma(a,(x,y,z))=\begin{cases}
 1,&\text{if } P_a(y,x)>P_a(z,y)\\
 0,&\text{if } P_a(y,x)<P_a(z,y)
 \end{cases}.$$
\end{definition}

In other words, $\Gamma(a,(x,y,z))$ is $0$ whenever the inversion of the pair $(y,x)$ occur before the inversion of the pair $(z,y)$ in the process of transforming the identity into the permutation $w$ by the action of each adjacent transposition corresponding to the letters of $a=a_1a_2\cdots a_{\ell(w)}$. Note also that since the subword $xyz$ of the identity is transformed into the subword $zyx$ of $w$ by the action of the letters of $a$, we must have $P_a(y,x)>P_a(z,x)>P_a(z,y)$ when $\Gamma(a,(x,y,z))=1$, and
$P_a(z,y)>P_a(z,x)>P_a(y,x)$ when $\Gamma(a,(x,y,z))=0$.
For instance, the analysis of the diagram in Example \ref{ex:action} for the reduced word $a=213213$ of the longest permutation of $\ss_4$ shows that $\Gamma(a,(123))=1$, $\Gamma(a,(124))=1$, $\Gamma(a,(134))=0$ and $\Gamma(a,(234))=0$.

In \cite{gms} it is shown that the map $\Gamma$ is invariant for the commutation classes of $C(\wo)$. We generalize this result for any permutation $w$.

\begin{proposition}\label{lemma:classinvariant}
 Two reduced words $a,b\in\Red(w)$ are in the same commutation class if and only if $\Gamma(a,(x,y,z))=\Gamma(b,(x,y,z))$, for all triple $(x,y,z)\in\tt_w$.
\end{proposition}
\begin{proof}
 If $a$ and $b$ are in the same commutation class, we may assume without loss of generality that there is an integer $i\in [\ell(w)-1]$ such that $\cc_{i}(P_a)=P_b$. The map $\cc_{i}$ changes the labels $i$  and $i+1$ of two cells that are not in the same row nor the same column. If one of these corresponds to a cell of a triple $(x,y,z)\in \tt_w$, then the change of the label does not modify the relative values of the end points of the cells $(z,y)$, $(z,x)$ and $(y,x)$. It follows that $\Gamma(a,(x,y,z))=\Gamma(b,(x,y,z))$ for all $(x,y,z)\in \tt_w$.

Reciprocally, suppose the $\Gamma$-value of $a$ and $b$ is the same for all triples in $\tt_w$.
Assume all labels  $1,\ldots,i-1$ are in the same cells in both tableaux $P_a$ and $P_b$, and $P_a(z,y)=i$ and $P_b(z,y)=i+k$ for some $k\geq 1$. Then, the permutation associated with $a_1\cdots a_{i-1}=b_1\cdots b_{i-1}$ has the factor $yz$.

If the label $i+k-1$ in $P_b$ is in row $z$, say in cell $(z,x)$, then the permutation associated with $b_1\cdots b_{i+k-1}$ has the factor $yxz$.  This means that $(y,x)$ is also an inversion for $w$, and therefore  $\Gamma(a,(x,y,z))=1$ and $\Gamma(b,(x,y,z))=0$, contradicting our assumption.
The same reasoning shows that the label $i+k-1$ cannot be in column $y$ of $P_b$, and therefore it must be in a cell which is not in row $z$ nor in column $x$. It follows that we can exchange the labels $i+k$ and $i+k-1$ using a commutation. That is $\cc_{i+k}(P_b)=P_{b'}$, such that all cells $1,\ldots,i$ are in same cells in both $P_a$ and $P_{b'}$.  Repeating the argument, there is a sequence of integers $i_1,\dots,i_j$ such that $\cc_{i_1}\cdots\cc_{i_j}(P_b)=P_a$, showing that $a$ and $b$ are in the same commutation class.
\end{proof}

\begin{definition}\label{def:distance}
 Given $a,b\in \Red(w)$, let

 $$ t(a,b)=\sum_{(x,y,z)\in \tt_w}\Gamma(a,(x,y,z))\oplus_2\Gamma(b,(x,y,z)),$$
 where $\oplus_2$ represents the sum  modulo 2.
\end{definition}

The number  $t(a,b)$ gives the number of triples in $\tt_w$ for which the $\Gamma$-value of $a$ and $b$ are distinct. Note that  by Proposition \ref{lemma:classinvariant}, we have $t(a,b)=t(a',b')$ for any $a'\in[a]$ and $b'\in [b]$. In particular, if $a$ and $b$ are in the same commutation class we get $t(a,b)=0$.

\begin{proposition}\label{lema:dist1}
 Let $a,b\in\Red(w)$. Then, $d([a],[b])=1$ if and only if $t([a],[b])=1$.
\end{proposition}
\begin{proof}
 If $d([a],[b])=1$, then we may assume without loss of generality that $a\lbr b$, which means that $\bb_i(P_a)=P_{b}$ for some integer $i$. It follows from Lemma \ref{lemma:admissiblehook} that the only triple having distinct $\Gamma$-values for $a$ and $b$ is $(x,y,z)$, $x<y<z$, corresponding to the cells $(z,y)$, $(z,x)$ and $(y,x)$ having labels $i-1$, $i$, and $i+1$, and therefore $t([a],[b])=1$.

 Assume now that $t([a],[b])=1$, and let $(x,y,z)\in \tt_w$ be the only triple having distinct $\Gamma$-values for $a$ and $b$. Using the same argument of the proof of Proposition \ref{lemma:classinvariant}, we may assume that all cells of $P_a$ and $P_b$ have the same value, with the exception of cells in positions $(z,y)$, $(z,x)$ and $(y,x)$. Suppose $T(a,(x,y,z))=0$ and $T(b,(x,y,z))=1$. Then, we must have $P_a(y,x)=i$, $P_a(z,x)=j$ and $P_a(z,y)=k$, with $i<j<k$, while $P_b(y,x)=k$, $P_b(z,x)=j$ and $P_b(z,y)=i$. We will show that if $j>i+1$ we can use commutative relations to  swap the integer $j$ with the integer $i+1$ in $P_a$.

 When  $j>i+1$  the permutation $a_1\cdots a_{i-1}$ has the factor $xy$ since $P_a(y,x)=i$, and $b_1\cdots b_{i-1}$ has the factor $yz$ since  $P_b(z,y)=i$. Since these two permutations coincide, we conclude that $xyz$ is a factor in this permutation, with $x<y<z$.
 If $j-1$ is in row $z$ of $P_a$, and also $P_b$, say in position $(z,u)$, then the permutation $a_1\cdots a_i\cdots a_{j-1}a_j$ implies the inversion of $(u,x)$ and $(u,y)$ by the action of $a_{i+1}\cdots a_{j-2}$. But then, we cannot have $P_{b}(z,u)=j-1$, contradicting our assumption. Thus, $j-1$ cannot be in row $z$ of $P_a$. The same argument shows that $j-1$ cannot be in column $x$ of $P_a$. Therefore, we can swap the integers $j$ and $j-1$ using the map $\bb_j$, corresponding to a commutation relation in $a$. Repeating the argument, we may assume that $j=i+1$, and by Lemma \ref{lemma:admissiblehook} we also have $k=i+2$.

 Therefore, we may use commutation relations to find reduced words $a',b'\in \Red(w)$ such that $a\sbr a'$ and $b\sbr b'$, where all cells of $P_{a'}$ and $P_{b'}$  have the same value, with the exception of cells in positions $(z,y)$, $(z,x)$ and $(y,x)$, where we have $P_{a'}(y,x)=i$, $P_{a'}(z,x)=i+1$, $P_{a'}(z,y)=i+2$, and $P_{b'}(y,x)=i+2$, $P_{a'}(z,x)=i+1$, $P_{a'}(z,y)=i$. It follows that $\bb_{i+1}(P_{a'})=P_{b'}$, that is $d([a],[b])=1$.
\end{proof}

 Proposition \ref{lema:dist1} shows that each triple in $\tt_w$ where two reduced decompositions $a,b\in\Red(w)$  have distinct $\Gamma$-values only for this triple, corresponds to a long braid relation between $a$ and $b$, or equivalently, to a mapping $\bb_i$ of $P_a$ into $P_{b}$. Also, it follows that $t([a],[b])\leq d([a],[b])$ for any $a,b\in\Red(w)$.

\begin{lemma}\label{lemma:min_max}
 The reduced word $a_{\min}$ (resp. $a_{\max}$) for $w\in\ss_n$ have $\Gamma$-value equal to $0$ (resp. $1$), for all triples in $\tt_w$.
\end{lemma}
\begin{proof}
We prove the result for $a_{\min}$ only. The other case is analogous.
Let $(x,y,z)$, $x<y<z$, be a triple in $\tt_w$ and consider cells $(y,x)$, $(z,x)$ and $(z,y)$  in the tableau $P_{a_{\min}}$. Since these cells correspond to inversions in $w$, the word $zyx$ is a subword of $w$. Thus, in the process of constructing $a_{\min}$, the descent corresponding to the pair $(z,y)$ will appear first, followed by the descent corresponding to $(z,x)$, and finally $(y,x)$. This means that the labels $i,j,k$ of cells  $(y,x)$,  $(z,x)$ and $(z,y)$ satisfy $i<j<k$, proving that $\Gamma(w,(x,y,z))=0$.
\end{proof}

\begin{definition}
 Let $w\in\ss_n$ and $a,b\in\Red(w)$. Denote by $\sup(a)$  the set of all triples $t\in \tt_w$ for which
 $\Gamma(a,t)=1$.
\end{definition}

By Lemma \ref{lemma:min_max}, we have $\sup(a_{\min})=\emptyset$ and $\sup(a_{\max})=\tt_w$.
That is, all braid relations in $a_{\min}$ have the form $i(i+1)i$ and all braid relations in $a_{\max}$ have the form $(i+1)i(i+1)$.
Moreover, we have $|\sup(a)|=t(a_{\min},a).$

\begin{lemma}\label{lemma:baixo}
 Let $a\in\Red(w)$ such that $a\notin [a_{\min}]$. Then, there is $b\in\Red(w)$ such that $[b]\lbr [a]$ and $\sup(b)\subset \sup(a)$.
\end{lemma}
\begin{proof}
 Assume there is no reduced word $b\in\Red(w)$ such that $[b]\lbr [a]$ and $\sup(b)\subset \sup(a)$. This means that $a$ cannot have a factor $(i+1)i(i+1)$. Therefore, each factor $(i+1)i$ in any word in the class $[a]$ can only be followed by a letter $i+1$ if there is a letter $i+2$ between them. But this implies $\Gamma(a, t)=0$ for all triple $t\in\tt_w$. That is, $a=a_{\min}$, contradicting our assumption.
\end{proof}

Note that in the condition of the lemma above, we have $\sup(b)=\sup(a)\setminus\{t\}$, for some $t\in\sup(a)$, and thus $d([a_{\min}],[b])=d([a_{\min}],[a])-1$.
Therefore, by successive applications of Lemma \ref{lemma:baixo} and Proposition \ref{lema:dist1} it follows that
\begin{equation}\label{eq:amin}
d([a_{\min}],[a])=t(a_{\min},a)=|\sup(a)|.
\end{equation}
An analogous result to Lemma  \ref{lemma:baixo} can be stated for the word $a_{\max}$, that is, if $a\notin [a_{\max}]$, then there is $b\in\Red(w)$ such that $[b]\lbr [a]$ and $\sup(a)\subset \sup(b)$. It follows that
 \begin{equation}\label{eq:amax}
 d([a],[a_{\max}])=t(a,a_{\max})=|\tt_w|-|\sup(a)|.
\end{equation}

This result shows that the map $t:\Red(w)\rightarrow [|\tt_w|]$ defined by  $t(a)=t(a_{\min},a)$
 is a rank function  for the graph $C(w)$, making it into a ranked partially ordered set with maximum and minimum. This partial order induces an orientation on $C(w)$.
\begin{proposition}
 Let $w\in\ss_n$. The partial order defined on the commutation classes of $C(w)$ given by the transitive closure of covering relations
 $$[a]<[b] \text{ if }a\lbr b\text{ and }t(b)=t(a)+1,$$
 makes $C(w)$ into a ranked partially ordered set with a unique minimal element $[a_{\min}]$ and a unique maximal element $[a_{\max}]$.
\end{proposition}

\section{Diameter of commutation graphs}

We can now give a formula for the diameter of $C(w)$, for any permutation $w\in\ss_n$.

\begin{theorem}\label{theorem:cardinality}
 The diameter of $C(w)$ is equal to the cardinality of $\tt_w$.
\end{theorem}
\begin{proof}
If $\tt_w$ is nonempty, then by
 Lemma \ref{lemma:min_max} we have $$d([a_{\min}],[a_{\max}])=|\tt_w|.$$

 Let $a,b\in\Red(w)$ such that $a\notin[b]$.
By equations \eqref{eq:amin} and \eqref{eq:amax}, we have
$$d([a_{\min}],[a])+d([a],[a_{\max}])+d([a_{\min}],[b])+d([b],[a_{\max}])=2|\tt_w|.$$ Using the triangle inequality, we conclude that
\begin{align*}
d([a],[b])&\leq \min\{d([a_{\min}],[a])+d([a_{\min}],[b]),\, d([a],[a_{\max}])+d([b],[a_{\max}])\}\\
&\leq |\tt_w|,
\end{align*}
proving that the distance between any two commutative classes $[a]$ and $[b]$ is at most $|\tt_w|$.
 Since this number is the largest possible distance between any two classes in $C(w)$, it is the diameter of the graph.
\end{proof}

\begin{remark}
The statistic used in the Section 4 induces a set-valued metric in the sense of \cite{reiner_roichman}. Proposition 3.12 in \cite{reiner_roichman} states that the diameter of a graph with a set-valued metric is equal to the maximal rank of the inclusion poset. Although we prove the same result for the graphs $C(w)$, in general $C(w)$ does not fulfill the requests of Proposition 3.12 in \cite{reiner_roichman} and it was necessary a more general proof.
 \end{remark}

 Let $w=w_1\cdots w_n\in\ss_n$ and let $p\in\ss_r$, for $r\leq n$. We say that $w$ {\it contains the pattern} $p$ if there exists a subsequence $w_{i_1}\cdots w_{i_r}$ whose elements are in the same relative order as the elements in $p$. If $w$ does not contain $p$, then we say that $w$ {\it avoids} $p$, or that $w$ is $p$-avoiding.

A permutation  having only one commutative class is said to be a fully commutative permutation. If $w\in\ss_n$ is fully commutative, the distance between any two reduced words for $w$ must be zero, which implies that $T_w$ is the empty set. That is, there is no triple $i<j<k$ with $w_i>w_j>w_k$, {\it i.e.} $w$ is $321$-avoiding. Thus we have recover a result of Billey, Jockusch, and Stanley \cite{321Avoiding}.

\begin{theorem}
 A permutation $w\in\ss_n$ is fully commutative if and only if it is $321$-avoiding.
\end{theorem}

A permutation  $w=w_1\cdots w_n$ is {\it unimodal} if there exists an index $i$, called the peak, such that $w_1<w_2<\cdots<w_i>w_{i+1}>\cdots>w_n$. Unimodal permutations  are characterized as avoiding the patterns $312$ and $213$, and are enumerated by $2^{n-1}$ \cite{ganon}. Using Theorem \ref{theorem:cardinality}, we can derive an explicit formula for the diameter of the commutation classes of unimodal permutations.

\begin{theorem}
 The diameter of the commutation graph of an unimodal permutation $w_1w_2\cdots w_n$ of $\ss_n$ with peak $i$ is
 $\ds\binom{n-i+1}{3}+\sum_{k=1}^{i-1}\binom{w_k-k}{2}$.
\end{theorem}
 \begin{proof}
If $w$ is unimodal with peak $i$, then the set $\tt_{w}$ is the union of sets $A\cup B$, where $A$  is formed by all triples $(w_{j_3},w_{j_2},w_{j_1})$ with $i\leq j_1<j_2<j_3\leq n$, and $B$ is formed by all triples $(w_k,w_{j_2},w_{j_3})$ with
$k<i<j_2<j_3\leq n$ and $w_k>w_{j_2}$. The set $A$ has precisely $\ds \binom{n-i+1}{3}$ elements, while in $B$ for each fixed $k$ there  are $\ds\binom{w_k-k}{2}$ triples, since there are precisely $w_k-k$ letters less than $w_k$ after the peak.
The result now follows from  Theorem \ref{theorem:cardinality}.
\end{proof}

 Since the longest permutation $\wo$ is the only unimodal permutation with peak $1$, we recover the following result from \cite{gms} (see also \cite{ziegler}) for the diameter of $C(\wo)$. This is the largest diameter for the commutation graph of a permutation in $\ss_n$.

 \begin{corollary}\label{cor:cardinality}
 The diameter of the commutation graph for the longest permutation $\wo$ of $\ss_n$ is $\ds\binom{n}{3}$.
\end{corollary}

The next proposition will help us to establish the maximal cardinalities that a graph $C(w)$ can have.

\begin{proposition}
  Let $\ds \ell:=\ell(\wo)=\binom{n}{2}$ be the length of the longest permutation of $\ss_n$, $\ds\delta:=\binom{n}{3}$ the diameter of $C(\wo)$, and  $w\in\ss_n$.
 If $\ell(w)=\ell -k$, then the diameter of $C(w)$ belongs to the interval $\ds \left[\delta-k(n-2),\delta-k(n-2)+\binom{k}{2}\right]$.
\end{proposition}

\begin{proof}
 If $\ell(w)=\ell -k$, then every pair $(a,b)$, with $1\leq a<b\leq n$ is transposed by $w$ except for $k$ pairs, $(a_i,b_i)$,
 with $1\leq i\leq k$. Thus, the set $\tt_w$ contains every triple $(a,b,c)$, with $1\leq a<b<c\leq n$, except for the ones for which two of the elements are $a_i$ and $b_i$. For each $i$ there are $n-2$ triples of this form, and then, by Theorem \ref{theorem:cardinality}, the diameter of $C(w)$ is at least $\delta -k(n-2)$.

 Each two pairs $(a_i,b_i)$ may have a common element, and if this is the case for every two pairs then the number of triples in $\tt_w$ is
 $\ds \delta - k(n-2)+\binom{k}{2}$, which imply that the diameter of $C(w)$ is at most $\ds \delta - k(n-2)+\binom{k}{2}$.
\end{proof}

When $n\geq 4$ the largest possible diameters for a graph $C(w)$, with $w\in\ss_n$, are, by decreasing order,
  $\delta, \delta -n+2, \delta -2n+5, \delta -2n+4$, corresponding to permutations with length $\ell$, $\ell-1$ and $\ell-2$. The longest permutation $\wo$ is the only permutation whose graph has length $\delta$.
  The unimodal permutations $(n-1)n(n-2)\cdots 21$ and $(n-2)(n-1)n(n-3)\cdots 21$,
  of lengths $\ell-1$ and $\ell-2$, are examples of permutations whose graphs have diameters $\delta -n+2$ and $\delta-2n+5$, respectively. Finally, the permutation $w=(n-1)n(n-2)(n-3)\cdots 4312$ has length $\ell-2$ and the graph $C(w)$ has diameter $\delta -2n+4$.

 %%%%%%%%%%%%%%%%%%%%%%%%%%%%%%%%%%%%%%%%%%%
 %%%%%%%%%%%%%%%%%%%%%%%%%%%%%%%%%%%%%%%%%%%
 %%%%%%%%%%%%%%%%%%%%%%%%%%%%%%%%%%%%%%%%%%%

\end{document}